\documentclass[11pt]{article}

\usepackage[margin=1in]{geometry}
\usepackage{amsmath,amsfonts,amssymb,amsthm}
\usepackage{graphicx}
\usepackage[authoryear]{natbib}
\usepackage[colorlinks,citecolor=blue,urlcolor=blue]{hyperref}

\newtheorem{lemma}{Lemma}

\newcommand{\gE}{\gamma_{\mathrm{E}}}
\newcommand{\gS}{\gamma_{\mathrm{score}}}

\begin{document}

\title{Rethinking Mean Square Error: Information, Generalized Estimation,
and the James--Stein Paradox}

\author{Paul Vos\thanks{Department of Public Health, East Carolina University,
Greenville, North Carolina. Email: \texttt{vosp@ecu.edu}. ORCID: 0000-0001-9996-5627.}}

\date{\today}

\maketitle

\begin{abstract}
The James--Stein estimator's dominance over maximum likelihood in mean
square error has been called a paradox because maximum likelihood is
known to be superior in many other respects. One response, due to
Efron, is to question maximum likelihood. Another is to question MSE.
We pursue the second and compare MSE with $\Lambda$-information
\citep{VosWu2025} as criteria for assessing estimators. The comparison
rests on two distinctions: between point estimators and generalized
estimators---functions of the sample and parameter jointly, with the
score as archetype---as inferential objects, and between pointwise and
family-aware assessment criteria. An elementary lemma shows that no
pointwise criterion, MSE or any other risk built from a loss function,
admits a uniformly optimal estimator; $\Lambda$-information, which is
family-aware and parameter-invariant, is uniformly maximized by the
score. A point estimator is assessed through the generalized
estimators it induces, and under the score map its
$\Lambda$-efficiency is the fraction of Fisher information the
statistic retains, placing the criterion in Fisher's information-loss
tradition. On unbiased estimators, $\Lambda$-efficiency coincides with
variance-based efficiency. Returning to James--Stein, the paradox
dissolves: maximum likelihood is fully efficient because it is
sufficient, while the James--Stein statistic is exactly two-to-one in
the sample, and the information it destroys---computed
exactly---is concentrated precisely where its MSE advantage is
greatest. MSE retains its proper domain under genuine squared-error
loss.
\end{abstract}

\section{Introduction}

In a recent article on James--Stein (JS) estimation, \citet{Efron2024-ef}
claims that minimizing mean square error (MSE) demonstrates ``in
an inarguable way, the virtues of shrinkage estimation'' and questions
whether ``two centuries of statistical theory, ANOVA, regression,
multivariate analysis, etc., depended on maximum likelihood estimation''
needs rethinking. Efron's response to the JS paradox---that JS dominates
the maximum likelihood (ML) estimator in MSE while ML remains superior
by many other criteria---is to question ML as the desired estimator.
Another response is to question MSE as the appropriate assessment
criterion. This paper pursues the second. We compare MSE with $\Lambda$-information
\citep{VosWu2025} as criteria for assessing estimators, and argue
that the JS paradox is more naturally read as a symptom of MSE's
limitations than as evidence against ML.

The comparison rests on a framework we develop in Section \ref{sec:Framework}:
a statistical model is a family $M$ of distributions whose smooth
manifold structure is independent of any particular parameterization,
and assessment criteria should respect this structure rather than
depend on the choice of labeling. With this framework in place, the
comparison between MSE and $\Lambda$-information unfolds along two
parallel lines. The first concerns inferential objects: $\Lambda$-information
evaluates generalized estimators---functions of the sample and the
parameter jointly, with the score as archetype---rather than point
estimators alone. Generalized estimators exist in cases where point estimators
do not, live in a Hilbert space that supports information-geometric
reasoning, and match Fisher's view of estimation as a continuum of
significance tests. A point estimator enters the comparison through
the generalized estimators it induces; under the score map, its
$\Lambda$-efficiency is the fraction of the Fisher information the
statistic retains, placing the criterion in the information-loss
tradition of \citet{Fisher1922-ie}. The second concerns assessment criteria themselves.
On the class of unbiased estimators, $\Lambda$-efficiency coincides
with variance-based efficiency; the two criteria agree where classical
theory is strongest. For biased estimators, $\Lambda$-information
remains parameter-invariant and uniformly bounded by Fisher
information, and the score attains the bound uniformly. MSE is
parameter-dependent, and an elementary impossibility lemma shows that
no pointwise criterion---MSE, or any other risk built from a loss
function---admits a uniformly optimal estimator; family-awareness,
not a better loss, is what uniform optimality requires.

Our argument is that $\Lambda$-information is better suited than MSE
for the broad class of parametric inference problems concerned with
identifying distributions within a family. MSE retains a proper domain
under genuine squared-error loss, where MSE is not a proxy for something
else but the loss itself, and moments continue to play their established
role in guiding the choice of $M$. The framework we propose concerns
how to assess inference procedures once $M$ has been selected; it
does not displace moment-based reasoning where moments are the right
tool.

In this telling, JS recedes to its proper role. JS raises a natural
question---\emph{if JS beats ML in MSE, does that mean we should use
JS?}---whose answer depends on which criterion is used. The general
comparison of MSE and $\Lambda$-information proceeds without JS. At
the end, JS returns to show that under $\Lambda$-information the paradox
dissolves: ML is fully efficient and the framework requires no change
to classical theory. The dissolution is exact rather
than rhetorical: maximum likelihood retains the full Fisher information
because it is sufficient, while the JS statistic is two-to-one in the
sample, and the information it discards is concentrated precisely
where its MSE advantage is greatest.

The paper proceeds as follows. Section \ref{sec:Framework} develops
the framework. Section \ref{sec:GE-vs-PE} compares generalized estimators
with point estimators as inferential objects. Section \ref{sec:Lambda-vs-MSE}
compares $\Lambda$-information with MSE as assessment criteria,
establishing an impossibility lemma for pointwise criteria. Section
\ref{sec:JS-return} returns to the JS example: the JS statistic is exactly two-to-one,
its information deficit is computed exactly, and the paradox
dissolves.

\section{Framework\label{sec:Framework}}

\subsection{Statistical models as smooth families}

A statistical model is a family $M$ of probability distributions sharing
a common support $\mathcal{Y}$ and a common $\sigma$-field $\mathcal{B}$.
For countable $\mathcal{Y}$, each $m\in M$ assigns probabilities to
points; for uncountable $\mathcal{Y}$, each $m$ has a density with
respect to a dominating measure $\mu$. It is convenient to treat a
distribution as a single object that can be written either as a density,
a function on $\mathcal{Y}$, or as a measure, a function on $\mathcal{B}$;
the two contain the same information and we write $m$ for both, context
making clear which is meant. In either case, all distributions in
$M$ share the same support, and the family is the set of points so
described.

For families used in parametric inference, $M$ carries additional
structure that makes it a smooth manifold.\footnote{We work throughout
with $M$ as an open smooth manifold admitting a global parameterization,
which lets us avoid the notion of coordinate charts. Boundary cases---such
as $p=0$ or $p=1$ in the binomial---are handled by restricting attention
to the interior of the parameter space, which avoids the technical
complications of manifolds with boundary.} The smooth structure consists
of relationships among nearby distributions---how one distribution
can be perturbed into another---described locally by the tangent spaces
$T_{m}M$. A \emph{parameterization} is a smooth map $\theta_{M}:M\to\Theta$
with smooth inverse, where $\Theta$ is an open subset of $\mathbb{R}^{k}$;
that is, $\theta_{M}$ is a diffeomorphism between $M$ and its image.
We distinguish the parameterization $\theta_{M}$, which is a function,
from a \emph{parameter} $\theta=\theta_{M}(m)\in\Theta$, which is
the label the parameterization assigns to a distribution $m$. The
role of the parameterization is to describe the smooth structure of
$M$ in tractable coordinates, so that calculations (derivatives,
integrals, optimization) can be performed using familiar tools on
$\Theta$.

The crucial point is that the smooth structure is a property of $M$
itself, not of any particular parameterization. Two parameterizations
$\theta_{M}$ and $\xi_{M}$ provide different coordinate systems for
the same manifold; the manifold structure is what they have in common.

\subsection{Attributes and parameterizations}

Two kinds of objects on $M$ should be distinguished. An \emph{attribute}
is a property of a single distribution. It can be defined without
reference to any family: the mean of a distribution, its variance,
its entropy, its mode, and its support are attributes of any distribution
that has them, whether or not the distribution sits in a family. Attributes
describe what individual distributions are like.

A \emph{parameterization}, as defined above,
is a diffeomorphism $\theta_{M}:M\to\Theta$ encoding the smooth structure
of $M$. It is not a property of any single distribution but a description
of how distributions in $M$ relate to one another, and any diffeomorphism
will serve.

The two kinds of object are easily conflated, because some parameterizations
happen to coincide with attributes. The mean parameterization of a
normal-mean family, for instance, assigns to each distribution its
own mean. But the coincidence is incidental: the mean parameterization
is a parameterization by virtue of being a diffeomorphism that encodes
the smooth structure, not by virtue of naming an attribute. Reparameterizing
by the log mean, or by any other diffeomorphism, changes the labels
without changing the structure described.

The Kullback--Leibler
divergence $\textnormal{KL}(m_{1},m_{2})$ takes two distributions
and returns a number describing how distinguishable they are. Like
an attribute, KL does not require that its distributions belong to
a family but, when $m_{1}$ and $m_{2}$ belong to a family $M$, it
is related to Fisher information, a characteristic of $M$.

Fisher information
describes how rapidly the distributions in $M$ separate from one
another near a point, and in this sense it is a characteristic of
the family rather than a property of any single distribution. The connection is local:
when $m_{1}$ and $m_{2}$ are nearby distributions in $M$ with parameters
$\theta_{1}=\theta_{M}(m_{1})$ and $\theta_{2}=\theta_{M}(m_{2})$,
\[
\textnormal{KL}(m_{1},m_{2})\approx\tfrac{1}{2}(\theta_{1}-\theta_{2})^{t}\,I(\theta_{1})\,(\theta_{1}-\theta_{2}),
\]
where $I$ is the Fisher information expressed in the $\theta_{M}$
parameterization. 

This distinction---between a single-distribution attribute and a characteristic
of the family---bears directly on the comparison to come. The Cram\'er--Rao
lower bound relates the variance of an estimator, which is an attribute
of a single distribution, to the Fisher information, which is a characteristic
of the family. The bound thereby conflates two objects of different
kinds. Keeping them separate is part of what motivates assessing estimators
by information rather than by variance or MSE. KL and the Fisher metric
reappear in Section \ref{sec:Lambda-vs-MSE}.

\subsection{The units-of-measurement analogy}

The role of parameterization on $M$ has a natural analogue in the
role of measurement units on the sample space.

Consider a sample space $\mathcal{Y}$ of waist circumferences. The
space carries structure---length, in particular---that is invariant
under the group of affine transformations: $y\mapsto ay+b$. To compute
with length we must choose a unit, but the choice of inches or centimeters
has no bearing on what family $M$ we select to model the population,
because the underlying structure (length) is independent of the labeling.
An inference procedure that gave different answers depending on whether
the data were recorded in inches or in centimeters would make the
choice of unit part of the inference problem---treating a feature
of the recording as though it were a feature of the population.

The same holds for $M$ and parameterizations. The smooth structure
of $M$ is invariant under the group of diffeomorphisms---under reparameterization.
To compute with the smooth structure we must choose a parameterization,
but the choice of $\theta_{M}$ versus $\xi_{M}=f\circ\theta_{M}$,
where $f$ is a diffeomorphism, has no bearing on which distribution
in $M$ we identify, because the smooth structure is independent of
the labeling. An inference procedure that gives different answers
depending on the choice of parameterization makes a particular parameter
part of the inference problem---treating a feature of the labeling
as though it were a feature of the family.

The two structures are different. Length is invariant under the affine
group; smoothness is invariant under the diffeomorphism group. The
diffeomorphism group is much larger than the affine group, which reflects
the fact that the smooth structure of $M$ is a richer object than
the metric structure of $\mathcal{Y}$. But the roles played by units
and parameterizations are parallel: in both cases, the structure of
interest is invariant, and the labels we use to describe it are not.

\subsection{Implications for assessment}

This framework has a direct consequence for the assessment of inference
procedures. An assessment criterion that depends on $M$ and the data---the
actual objects of inference---but not on the choice of labeling is
one that respects the structure of the problem. A criterion that depends
on the choice of parameterization is one that is sensitive to the
labels rather than the structure; a criterion that ignores $M$ itself,
assessing a procedure one distribution at a time as though the rest of
the family were absent, discards the structure rather than respecting
it. The first failing is a lack of parameter invariance; the second is
a lack of family-awareness. Section \ref{sec:Lambda-vs-MSE} shows that
MSE exhibits both.

We take this as a desideratum; the two sections that follow develop
its consequences.

The parallel between the two informations is worth stating at the
outset. Just as Fisher information describes a property of the family
$M$---how rapidly its distributions separate near a point---$\Lambda$-information
describes a property of a generalized estimator, an object that shows
how the data are related to all the distributions in $M$. Both are
parameter-invariant: their values attach to the objects they describe
(the family, the estimator) and not to the parameterization used to
compute them.

\section{Generalized estimators and point estimators\label{sec:GE-vs-PE}}

The first of our two comparisons concerns the objects of inference
themselves. We first fix terminology. An \emph{estimator} is a function
on the sample space $\mathcal{Y}$; an \emph{estimate} is its value
at a particular sample $y$. A point estimator is thus a random variable,
and a point estimate $\hat{\theta}(y)$ is a real number---the single
distribution in $M$ labeled $\hat{\theta}(y)$. A \emph{generalized
estimate} for $y$ is a function on the parameter space $\Theta$,
so that a \emph{generalized estimator} is a function on $\mathcal{Y}\times\Theta$.

A point estimator induces a generalized estimator, and it does so in
two ways that must be distinguished. The first map uses expectation:
\begin{equation}
\gE(\hat{\theta})=g_{E(\hat{\theta})}=\hat{\theta}(y)-E_{\theta}(\hat{\theta}),\label{eq:gammaE}
\end{equation}
which depends on $\theta$ through the expectation. The second map uses
the model for the estimator itself. The statistic $\hat{\theta}$ has,
at each $m\in M$, a marginal distribution; these distributions form a
statistical model in their own right, with log-likelihood
$\ell_{\hat{\theta}}(y,\theta)$, the log of the density of
$\hat{\theta}$ evaluated at $\hat{\theta}(y)$. Its score defines the
second map:
\begin{equation}
\gS(\hat{\theta})=g_{\hat{\theta}}=\ell'_{\hat{\theta}}.\label{eq:gammaS}
\end{equation}
For a fixed sample $y$, the point estimate $\hat{\theta}(y)$ is one
value, while $g_{E(\hat{\theta})}(y,\cdot)$ and
$g_{\hat{\theta}}(y,\cdot)$ take a value at every point of $\Theta$. The notation $g$ follows \citet{VosWu2025}, who
introduce generalized estimators in general form.\footnote{In
\citet{VosWu2025} the symbol $g_{\hat{\theta}}$
denotes $\hat{\theta}-E_{\theta}(\hat{\theta})$, the only induced
generalized estimator considered there. The present paper refines
that notation: with two maps in play, the expectation-map image
carries the subscript $E(\hat{\theta})$, and $g_{\hat{\theta}}$ is
reserved for the image under the score map.} The maps and their
images must be kept distinct: $\gE$ and $\gS$ are functions whose
common domain is the point estimators, while their images
$g_{E(\hat{\theta})}$ and $g_{\hat{\theta}}$ are generalized
estimators---functions on $\mathcal{Y}\times\Theta$.

The two maps behave differently, and the differences matter for what
follows. First, $\gS$ depends on $\hat{\theta}$ only through the
$\sigma$-field the statistic generates: if $\tau$ is injective, the
density of $\tau(\hat{\theta})$ is the density of $\hat{\theta}$
multiplied by a Jacobian free of $\theta$, so
\[
g_{\tau(\hat{\theta})}=\gS(\tau(\hat{\theta}))=\gS(\hat{\theta})=g_{\hat{\theta}}.
\]
The expectation map has no such invariance: for nonlinear $\tau$,
$\tau(\hat{\theta})-E_{\theta}\tau(\hat{\theta})$ is not a rescaling of
$\hat{\theta}-E_{\theta}(\hat{\theta})$, and the two images
$g_{E(\tau(\hat{\theta}))}$ and $g_{E(\hat{\theta})}$ are different
generalized estimators. Second, Section \ref{sec:Lambda-vs-MSE} shows
that the assessment of $g_{\hat{\theta}}$ is exactly the Fisher
information carried by the statistic $\hat{\theta}$, connecting
generalized estimation to the information-loss program of
\citet{Fisher1922-ie} and its development in curved exponential
families by \citet{Efron1975-ic}. Third, the expectation map remains
the operational choice when only the moment functions
$E_{\theta}(\hat{\theta})$ and $V_{\theta}(\hat{\theta})$ of the
estimator are available.

The archetype of a generalized estimator is not built from a point
estimator at all. For scalar parameter $\theta$, the score $\ell'=d\ell(y,\theta)/d\theta$
is a function on $\Theta$ for each sample, and it is the generalized
estimator against which others are measured. Write $Y$ for the sample
regarded as a random variable. Under mild regularity conditions the
score $\ell'(Y)$ has mean zero and finite variance at each $\theta$,
so it is square-integrable; evaluated at $\theta$ it is then an element
of a Hilbert space $H_{\theta}M$ whose subspace spanned by the score
components is the tangent space $T_{\theta}M$. For a scalar parameter,
$\ell'(Y)$ at $\theta$ is a single such vector; for a vector parameter,
the components of the gradient $\nabla\ell(Y)$ span $T_{\theta}M$.
The notation is parameter-indexed, but the score and the quantities
built from it are parameter-invariant: a reparameterization relabels
the same elements of $H_{\theta}M$. In point estimation the score
describes how the maximum likelihood estimate changes near a distribution;
as a generalized estimator, the score \emph{is} the estimator. In
generalized estimation the maximum likelihood estimate indicates where
the score crosses the horizontal axis. 
In a full exponential family the maximum likelihood estimator is an
injective function of the sufficient statistic, so under the score map
$g_{\hat{\theta}}=\gS(\hat{\theta})$ is the full-data score in
\emph{every} parameterization.

Generalized estimators are formally close to the estimating functions
of \citet{Godambe1974-pt} (see also \citealp{Godambe1978-jw,Godambe1989-yo}):
both are functions on $\mathcal{Y}\times\Theta$ satisfying
$E_{\theta}(g)=0$, and in both theories the score is the optimal member
of the class. The difference lies in what the function is for. In
estimating-equation theory the function is scaffolding: its role is to
produce a point estimator as the root of $g(y,\theta)=0$, and the
optimality of $g$ is valued because it transfers to that root. In
generalized estimation the function itself is the inferential object.
There is no obligation to solve $g=0$; the value $g(y,\theta_{\circ})$
is read directly as evidence discriminating $m_{\theta_{\circ}}$ from
its neighbors in $M$. A generalized estimator does not so much
estimate a parameter as discriminate among the distributions of $M$,
and this shift is what allows samples at which $g$ has no root---the
boundary samples of Section \ref{subsec:Existence}---to remain fully
interpretable.

Generalized estimators are a richer class than point estimators, and
the richness is not idle. We note three respects in which generalized
estimators behave better as objects of inference.

\subsection{Existence\label{subsec:Existence}}

A point estimator must return a distribution in $M$. For the maximum
likelihood estimator this means returning the maximizer of the likelihood
over $\Theta$, and when that maximizer lies outside $\Theta$, no estimate
exists. Consider the binomial family with support $\{0,1,\ldots,n\}$
and parameter space $\Theta=(0,1)$. At $y=0$ the likelihood $p\mapsto(1-p)^{n}$
is maximized at $p=0$, and at $y=n$ the likelihood $p\mapsto p^{n}$
is maximized at $p=1$; neither maximizer lies in the open set $\Theta$.
The maximum likelihood point estimate fails to exist at these samples.

The score function does not fail. For the binomial, $\ell'(y,\theta)$
is defined for every $y\in\{0,1,\ldots,n\}$ and every $\theta\in\Theta$,
including $y=0$ and $y=n$. The generalized estimator exists where
the point estimator does not, and inference can proceed.

This is not a peculiarity of the binomial. Whenever a sample lies
at the edge of what the family can accommodate, point estimation breaks
down while the score remains available. The generalized estimator's
existence is tied to the smooth structure of $M$---the differentiability
of the log-likelihood---rather than to the solvability of an optimization
problem on $\Theta$.

A natural repair suggests itself: close the parameter space by
appending $p=0$ and $p=1$ to $\Theta$, so that the boundary samples
have maximum likelihood estimates. The repair has costs on both
sides. On the theory side, $M$
ceases to be an open manifold of mutually absolutely continuous
distributions: the appended points are degenerate, the common-support
structure of Section \ref{sec:Framework} fails, and results split into
interior cases with boundary exceptions where the open family needs a
single statement. On the applied side, the appended points often
describe no population under study: in many applications it is known
that the population does not consist entirely of successes or entirely
of failures, so $p=0$ and $p=1$ are not distributions the model
intends. Closing the space is then a crutch for an assessment theory
that requires a point estimate to exist, not a statement about the
phenomenon being modeled. The generalized estimator requires no such
repair.

\subsection{Linear structure}

At each parameter value $\theta$, there is a Hilbert space that contains
the tangent space of $M$ as well as all the generalized estimators
$g$ evaluated at $\theta$. The Hilbert space inner product supports
orthogonal decomposition, which offers a route to nuisance parameters:
estimators can be orthogonalized within the space rather than by seeking
an orthogonal parameterization of $\Theta$, which may not exist. More
broadly, the Hilbert space is what makes the information-geometric
apparatus available: angles, projections, and the Pythagorean relationships
that underlie the additivity of information. Outside of the estimation
of location parameters for normally distributed data, this linear
structure has no counterpart for point estimators.

The geometric perspective has deep roots. \citet{rao1945information}
first recognized that families of distributions form differentiable
manifolds with Fisher information as a Riemannian metric. \citet{Efron1975-ic}
characterized statistical models through their curvature. \citet{Amari1990-bo}
developed the comprehensive differential-geometric framework that
established information geometry as a field. \citet{VosWu2025} and
\citet{Vos2025geometry} use the Hilbert bundle over $M$ to study
generalized estimators.

\subsection{Unification of estimation and testing\label{subsec:unification}}

A generalized estimator is, at each fixed parameter value $\theta_{\circ}$,
a test statistic for $H_{\circ}:\theta=\theta_{\circ}$. The function
$g(\cdot,\theta_{\circ})$ takes a value at every point
of $\mathcal{Y}$, and its magnitude orders the sample space: larger
values of $|g(y,\theta_{\circ})|$ indicate the sample
is farther in the tails when $\mathcal{Y}$ is ordered by $|g(y,\theta_{\circ})|$.
A generalized estimator is thus a smooth continuum of test statistics,
one for each candidate distribution.

This matches Fisher's view of estimation as a continuum of significance
tests. A point estimator answers the question ``which distribution?'';
a generalized estimator answers the family of questions ``how well
does each distribution account for this sample?'' The two are not
in tension---the point estimate is recovered as the distribution at
which the generalized estimator's test statistic is least extreme---but
the generalized estimator carries more. Estimation and testing, often
treated as separate activities, are two readings of the same object.

\section{\texorpdfstring{$\Lambda$-Information and mean square error}{Lambda-Information and mean square error}\label{sec:Lambda-vs-MSE}}

The second comparison concerns the criteria by which estimators are
assessed. We define $\Lambda$-information, record its agreement with
variance-based efficiency on unbiased estimators, and establish an
impossibility lemma showing that no pointwise criterion---MSE among
them---admits a uniformly optimal estimator. The comparison with MSE
then proceeds through the properties the lemma organizes: uniform
optimality, family dependence, and parameter invariance, together
with a control case showing that invariance alone is not enough. The James--Stein example is deferred to
Section \ref{sec:JS-return}; here the comparison rests on the criteria
themselves.

Let $g$ be a generalized estimator for scalar parameter $\theta$---a
function on $\mathcal{Y}\times\Theta$ with $E_{\theta}(g)=0$ for all
$\theta$, square-integrable at each $\theta$, and smooth in $\theta$;
these are the conditions of \citet[Definition 1]{VosWu2025}, and every
generalized estimator below is assumed to satisfy them.
The $\Lambda$-information of $g$ is the square of the mean rate at
which $g$ changes with the parameter, standardized by its variance:
\begin{equation}
\Lambda(g)=\frac{\left(E(g')\right)^{2}}{V(g)}.\label{eq:Lambda_g}
\end{equation}
The standardization is necessary because the raw rate of mean change
can be made arbitrarily large by rescaling $g$; only the standardized
rate is meaningful.\footnote{A note on notation: the expectation
operator no longer carries the subscript $\theta$ because it, along
with $g$, the variance $V(g)$, $\Lambda$, and Fisher information $I$,
are all functions of a parameter. Notation displaying a particular
parameter would draw attention to a choice that is
immaterial---all of these can be viewed as functions on $M$. The
Kullback--Leibler divergence of Section \ref{subsec:MKL} again carries
a subscript, $E_{m_{1}}$, because it has two arguments and the
notation must indicate which supplies the expectation.}

A point estimator $\hat{\theta}$ is assessed through the generalized
estimators it induces; Section \ref{sec:GE-vs-PE} distinguished the
score map from the expectation map. For the expectation-map image
$g_{E(\hat{\theta})}$, (\ref{eq:Lambda_g}) becomes
\begin{equation}
\Lambda(g_{E(\hat{\theta})})=\frac{\left((E\hat{\theta})'\right)^{2}}{V(\hat{\theta})}.\label{eq:Lambda_theta_hat}
\end{equation}
For a point estimator so assessed, the square root of the
$\Lambda$-information measures how fast the standardized mean of the
estimator changes with the parameter.

$\Lambda$-information is bounded above by Fisher information, $\Lambda(g)\le I$,
with equality for the score. This follows from differentiating the
identity $E(g)=0$ in $\theta$, which gives $E(g')+E(g\,\ell')=0$,
hence $E(g')=-\mbox{Cov}(g,\ell')$; the Cauchy--Schwarz
inequality then gives the bound. (The relation $E(g')=-\mbox{Cov}(g,\ell')$
is an instance of the second Bartlett identity \citep{Bartlett1953},
which for the score itself reads $E(\ell'')+E\{(\ell')^{2}\}=0$.) The
$\Lambda$-efficiency of $g$ is
\begin{equation}
\mbox{Eff}^{\Lambda}(g)=\frac{\Lambda(g)}{I}\le1.\label{eq:Eff_g}
\end{equation}
The efficiency and its multivariate form are detailed in \citet{VosWu2025};
we recall in Section \ref{subsec:bound} only what the James--Stein
example requires.

The score map gives $\Lambda$ a second reading, as exact information
accounting. For the score-map image $g_{\hat{\theta}}$, the identities above applied to
the marginal model of $\hat{\theta}$ give
$E(g')=-I_{\hat{\theta}}$ and $V(g)=I_{\hat{\theta}}$, where
$I_{\hat{\theta}}$ is the Fisher information of the marginal model of
the statistic, so
\begin{equation}
\Lambda(g_{\hat{\theta}})=I_{\hat{\theta}},\qquad
\mbox{Eff}^{\Lambda}(g_{\hat{\theta}})=\frac{I_{\hat{\theta}}}{I}\le1.\label{eq:Lambda_score_map}
\end{equation}
Under the score map, the $\Lambda$-efficiency of a point estimator
is the fraction of the Fisher information that the statistic retains,
with equality exactly when $\hat{\theta}$ is sufficient. This is
Fisher's measure of the information lost in reducing the sample to a
statistic \citep{Fisher1922-ie}, computed by \citet{Efron1975-ic} for
maximum likelihood in curved exponential families; the bound
(\ref{eq:Eff_g}) becomes the statement that a statistic cannot carry
more information than the sample. With the definitions in hand, the
comparison with MSE proceeds through the subsections that follow.

\subsection{Agreement on unbiased estimators}

Unbiasedness constrains the mean function of the estimator, which
the score-map image ignores entirely, so this subsection concerns the
expectation-map image.
For an unbiased estimator, $E(\hat{\theta})=\theta$, so
$(E\hat{\theta})'=1$ and $\Lambda(g_{E(\hat{\theta})})=1/V(\hat{\theta})$.
The $\Lambda$-efficiency is then $1/\{I\cdot V(\hat{\theta})\}$, which
is exactly the classical efficiency defined through variance. On the
class of unbiased estimators, $\Lambda$-information and variance-based
assessment coincide.

The information bound here recovers the Cram\'er--Rao inequality,
but in a form that keeps its two sides distinct. Written as $\Lambda(g_{E(\hat{\theta})})\le I$,
the left side assesses the generalized estimator $g_{E(\hat{\theta})}$,
an object that depends on $M$, while the right side is the Fisher
information, a characteristic of $M$. The classical statement $V(\hat{\theta})\ge1/I$
instead places the variance of $\hat{\theta}$---an attribute of a
single distribution---on one side and a family characteristic on the
other, joining two objects of different kinds. The $\Lambda$ form
relates like to like: an assessment of an $M$-dependent estimator,
bounded by a characteristic of $M$.

The two criteria agree where classical theory is strongest, and $\Lambda$-information
departs from variance only for biased estimators---where variance-based
efficiency is silent, because the Cram\'er--Rao bound depends on the
bias function and so compares an estimator only against others sharing
its bias.

\subsection{Pointwise criteria and an impossibility lemma\label{subsec:impossibility}}

The comparison in this subsection is stated on $M$ rather than on
$\Theta$: a point estimator $\hat{\theta}$ determines the estimator
$\hat{m}=\theta_{M}^{-1}(\hat{\theta})$, the distribution labeled by
the point estimate. The move costs nothing and buys generality, for
none of the smooth structure of Section \ref{sec:Framework} is used
here: $M$ need only be a set of distributions with common support, and
$\theta_{M}$ need only be an injection---a labeling---rather than a
parameterization. Call an assessment criterion \emph{pointwise} if it
assigns to an
estimator $\hat{m}$ and a single distribution $m\in M$ a risk
$R(\hat{m},m)\ge0$ computed from the sampling distribution of $\hat{m}$
under $m$ alone, with $R(\hat{m},m)=0$ if and only if $\hat{m}=m$
almost surely under $m$. MSE is pointwise:
$E_{m}||\hat{\theta}-\theta_{M}(m)||^{2}$ uses $m$ and the distribution of
$\hat{\theta}$ at $m$, and nothing else. So is any risk built from a
loss function $L$ with $L(m,m)=0$ and $L>0$ off the
diagonal---standardized quadratic loss, divergence-based losses, and
the mean Kullback--Leibler divergence of Section \ref{subsec:MKL}
among them.

\begin{lemma}\label{lem:pointwise}
Let $M$ be a set of distributions with common support, containing more
than one member, and let $R$ be pointwise. Then no estimator is
uniformly $R$-optimal over $M$.
\end{lemma}

\begin{proof}
Fix $m_{0}\in M$. The constant estimator $\hat{m}\equiv m_{0}$ has
$R(m_{0},m_{0})=0$, so a uniformly optimal $\hat{m}^{*}$ must satisfy
$R(\hat{m}^{*},m_{0})=0$, hence $\hat{m}^{*}=m_{0}$ almost surely
under $m_{0}$. The same argument applies at any $m_{1}\neq m_{0}$.
Since $m_{0}$ and $m_{1}$ share support, they are mutually absolutely
continuous, so $\hat{m}^{*}$ would be almost surely equal to two
different constants on the same sample space---a contradiction.
\end{proof}

The lemma is elementary, and that is its point. It shows that the
absence of a uniformly MSE-optimal estimator is not a defect of
squared error that a better loss might repair: no reweighting,
standardization, reparameterization-invariant loss, or divergence
escapes it, because all remain pointwise. Within the pointwise class
the standard responses change the question---restrict to unbiased
estimators, or replace uniform optimality with admissibility or
minimaxity---each importing considerations beyond the family and the
data. A criterion that escapes the lemma must use, at $m$, information
beyond the sampling distribution at $m$; it must be
\emph{family-aware}, responding to how the estimator behaves as $m$
moves within $M$.

\subsection{Uniform optimality}

$\Lambda$-information escapes Lemma \ref{lem:pointwise} because it is
not pointwise. Its numerator $E(g')$ is the response of $g$
to movement of the distribution within $M$; it is undefined at an
isolated distribution and uses precisely the structure the lemma shows
is needed. The escape can be seen at the lemma's own witness: the
constant estimator, which achieves zero risk somewhere under any
pointwise criterion, has mean response zero everywhere and hence
$\Lambda=0$---the worst value, not the best.

The bound $\Lambda(g)\le I$ holds at every parameter value, and the
score attains it at every parameter value. Optimality under
$\Lambda$-information is therefore uniform: the optimal generalized
estimator is the same---the score---across the entire parameter space.
No auxiliary concepts are required. Where pointwise criteria must
bring in unbiasedness, admissibility, or minimaxity to break ties the
criterion itself cannot break, information-based assessment orders the
class of generalized estimators against a single uniform benchmark.
That the benchmark is the score is itself Cauchy--Schwarz---$\mbox{Eff}^{\Lambda}(g)$
is the squared correlation between $g$ and the score
\citep[Corollary 1]{VosWu2025}; the content lies in the possibility,
which Lemma \ref{lem:pointwise} denies to every pointwise criterion.

\subsection{Family dependence}

MSE is defined at a single distribution: given a distribution and
an estimator, $E(\hat{\theta}-\theta)^{2}$ is computable without reference
to any family. $\Lambda$-information is not---it is built from the
score, which exists only relative to a family $M$, and from the rate
of change of the estimator's mean as the distribution moves within
$M$. $\Lambda$-information cannot be defined at an isolated distribution.

This dependence is a feature. Parametric inference is the problem
of identifying a distribution within a specified family; the family
is part of the problem, not incidental to it. A criterion that requires
the family is a criterion that uses all the structure the problem
provides. A criterion defined without the family---MSE---discards
that structure, assessing the estimator as though the inferential
task were to locate a point in isolation rather than to distinguish
among the members of $M$.

\subsection{Parameter invariance\label{subsec:invariance}}

$\Lambda$-information attaches to the generalized estimator as an
object on $M$, not to the coordinates used to compute it. The precise
statement, in the multivariate form needed later, is the following;
part (i) also settles the sense in which the assessment ignores
$\theta$-dependent rescaling.

\begin{lemma}\label{lem:invariance}
Let $g$ be a generalized estimator with values in $\mathbb{R}^{k}$,
and write
$\Lambda(g)=C^{{\tt t}}V^{-1}(g)\,C$ with
$C_{ij}=\mbox{Cov}(g_{i},\partial_{j}\ell)=-E(\partial g_{i}/\partial\theta_{j})$,
as in (\ref{eq:Lambda_mu}) below. Then:
(i) for any invertible matrix function $A(\theta)$,
$\Lambda(A(\theta)\,g)=\Lambda(g)$;
(ii) under a reparameterization $\xi=\tau(\theta)$ with Jacobian
$J=\partial\theta/\partial\xi$, $\Lambda$ and the Fisher information
transform identically, $\Lambda_{\xi}=J^{{\tt t}}\Lambda_{\theta}J$
and $I_{\xi}=J^{{\tt t}}I_{\theta}J$, so the eigenvalues of
$\mbox{Eff}^{\Lambda}(g)=I^{-1/2}\Lambda(g)I^{-1/2}$---in particular
the scalar efficiency---do not depend on the parameterization.
\end{lemma}

\begin{proof}
(i) $g\mapsto Ag$ sends $C\mapsto AC$ and $V\mapsto AVA^{{\tt t}}$,
and $(AC)^{{\tt t}}(AVA^{{\tt t}})^{-1}(AC)=C^{{\tt t}}V^{-1}C$.
(ii) $g$ is unchanged as a function on $\mathcal{Y}\times M$ while
$\nabla_{\xi}\ell=J^{{\tt t}}\nabla_{\theta}\ell$, giving
$C_{\xi}=C_{\theta}J$ with $V$ unchanged, hence
$\Lambda_{\xi}=J^{{\tt t}}\Lambda_{\theta}J$; the same relation for
$I$ is standard. Then
$I_{\xi}^{-1}\Lambda_{\xi}=J^{-1}(I_{\theta}^{-1}\Lambda_{\theta})J$
is a similarity transformation, preserving eigenvalues.
\end{proof}

Consequently $\mbox{Eff}^{\Lambda}(g)$ is well defined for $g$ as an
object on $M$: if $g_{1}$ and $g_{2}$ are the same generalized
estimator expressed in the coordinates of $\theta_{M}$ and $\xi_{M}$
respectively, then
$\mbox{Eff}^{\Lambda}(g_{1})=\mbox{Eff}^{\Lambda}(g_{2})$. The
efficiency attaches to the object on $M$, not to the coordinate
labels---just as, in Section \ref{sec:Framework}, length attaches to a
waist and not to the choice of inches or centimeters.

Two different invariances are in play and must not be conflated. The
first attaches to a fixed generalized estimator as an object on $M$:
by Lemma \ref{lem:invariance}, the assessment of a given $g$---and in
particular of a given image $g_{E(\hat{\theta})}$---is the same in
every parameterization of $M$. The second concerns the maps, where the
question is invariance under relabeling the \emph{estimator}: for
injective $\tau$, do $\hat{\theta}$ and $\tau(\hat{\theta})$ receive
the same assessment? Under the score map they do,
$g_{\tau(\hat{\theta})}=g_{\hat{\theta}}$ (Section
\ref{sec:GE-vs-PE}), so the efficiency of a point estimator is a
property of the $\sigma$-field the statistic generates. Under the
expectation map they need not: for nonlinear $\tau$ the images
$g_{E(\tau(\hat{\theta}))}$ and $g_{E(\hat{\theta})}$ are different
generalized estimators, each with a parameter-invariant assessment of
its own, and the two assessments need not agree. The expectation map
thus makes the choice of scale for the estimator part of the
assessment, while leaving intact the invariance on $M$ of each image.

MSE is not parameter-invariant. The MSE of an estimator of $\theta$
differs from the MSE of the corresponding estimator of $\xi=\tau(\theta)$,
even though both name the same distribution, because squaring a
deviation is a coordinate-dependent operation and $\tau$ is generally
nonlinear.

\subsection{A control case: mean Kullback--Leibler divergence\label{subsec:MKL}}

Lemma \ref{lem:pointwise} asserts that no pointwise criterion, however
well chosen, admits a uniformly optimal estimator. The assertion is
best tested against the strongest available candidate: a criterion
that repairs MSE's remaining defects---its dependence on units and on
the parameterization---while staying within the pointwise class. This
subsection constructs that candidate and watches it fail, which
locates the property that carries the weight.

A practical objection to MSE, raised by Efron's own example, is that
total MSE can sum quantities with incommensurable units. Using Efron's
example, total MSE sums squared Daltons, squared meters, and
dimensionless counts. Ignoring units is not an option: a rescaling
that is immaterial for a single component---meters to
kilometers---changes the relative weighting of the components and so
can change which estimator total MSE prefers. Cost-based conversion to
common units (e.g., dollars squared) requires three arbitrary
conversion factors: Daltons, meters, and counts to dollars.

This unit problem arises because MSE assesses estimators via
distribution attributes, specifically means. One familiar remedy is to
standardize each quantity before summing, replacing the squared error
in component $i$ by $(\hat{\mu}_{i}-\mu_{i})^{2}/\sigma_{i}^{2}$; the
standardized terms are dimensionless and can be added. For the normal
family this standardized squared error is exactly twice the
Kullback--Leibler divergence, so the remedy coincides with the mean KL
divergence constructed here. KL divergence is the more general
construction: it requires no assumption of normality and is defined
for any family with common support. We therefore replace the mean
estimate $\hat{\mu}$ with the point in $M_{k}$ having that mean,
$\hat{m}=m_{\hat{\mu}}=\mu^{-1}(\hat{\mu})$, where $M_{k}$ is the
homoscedastic normal family of Section \ref{sec:JS-return}, equation
(\ref{eq:M}), and $\mu$ its mean parameterization. The two estimators
live in different spaces: $\hat{\mu}$ takes values in
$\mathbb{R}^{k}$, where points are compared by squared distance
$||\mu_{1}-\mu_{2}||^{2}$, while $\hat{m}$ takes values in $M_{k}$,
where points are compared by KL divergence:
\begin{align}
K\!L(m_{1},m_{2}) & =E_{m_{1}}(\log(m_{1}/m_{2}))\nonumber \\
 & =\frac{1}{2}||\mu(m_{1})-\mu(m_{2})||^{2}\label{eq:KL2nd}
\end{align}
The first line defines KL divergence generally; the second evaluates
it on $M_{k}$. Crucially, KL divergence has no units---logarithms are
dimensionless, and log differences are log ratios where units cancel.

Mean KL divergence (MKL) of $\hat{m}$ is:
\begin{align}
M\!K\!L(\hat{m}) & =E_{m}K\!L(\hat{m},m)\nonumber \\
 & =\frac{1}{2}E_{m}||\mu(\hat{m})-\mu(m)||^{2}\label{eq:MKL_MSE}
\end{align}
The first line defines MKL generally; the second holds for $M_{k}$.

Equation (\ref{eq:MKL_MSE}) shows that on $M_{k}$, MSE in the mean
parameterization equals twice MKL. MKL is therefore everything a
repaired MSE could hope to be: unit-free, defined directly on
distributions rather than on their labels, and invariant under
reparameterization of $M$. Yet it is pointwise. $K\!L(\hat{m},m)$ is
a functional of the pair $(\hat{m},m)$ alone; its arguments need not
belong to any family, and no structure of $M$ enters its computation.
Lemma \ref{lem:pointwise} applies, and, concretely, the James--Stein
dominance reported in Table \ref{tab:MSE-for-James-Stein} of Section
\ref{sec:JS-return} carries over to MKL verbatim, halved.

The control case sharpens the diagnosis of this section. Invariance
and unit-freedom are necessary for a criterion that respects the
structure of Section \ref{sec:Framework}, but they are not sufficient:
MKL has both and inherits MSE's failures anyway, because it, like MSE,
assesses a point estimator against the truth one distribution at a
time. The load-bearing property is family-awareness, and
family-awareness is not available to any criterion whose argument is a
point estimator assessed pairwise. It requires the change of
inferential object made in Section \ref{sec:GE-vs-PE}: from the point
estimator, which assigns to each sample a single distribution, to the
generalized estimator, which at each sample assigns a value to every
distribution in $M$ in order to discriminate among them.

\subsection{The multivariate information bound\label{subsec:bound}}

The James--Stein example of Section \ref{sec:JS-return} requires the
multivariate form of $\Lambda$-information; we record it here and
refer to \citet{VosWu2025} for derivations. For a generalized
estimator $g$ with values in $\mathbb{R}^{k}$, the information is the
$k\times k$ matrix
\begin{equation}
\Lambda(g)=C^{{\tt t}}V^{-1}(g)\,C,\qquad
C_{ij}=\mbox{Cov}(g_{i},\partial_{j}\ell)=-E\left(\partial g_{i}/\partial\theta_{j}\right),\label{eq:Lambda_mu}
\end{equation}
the form used in Lemma \ref{lem:invariance}. For the expectation-map
image $g_{E(\hat{\mu})}=\hat{\mu}-E(\hat{\mu})$ of an estimator of
$\mu\in\mathbb{R}^{k}$, $C_{ij}=\mbox{Cov}(\hat{\mu}_{i},\partial_{j}\ell)=\partial E(\hat{\mu}_{i})/\partial\mu_{j}$
and $V(g_{E(\hat{\mu})})=V(\hat{\mu})$, giving the multivariate form
of (\ref{eq:Lambda_theta_hat}).
The scalar information $\lambda(g)=\mbox{tr}\,\Lambda(g)$
sums the eigenvalues of $\Lambda(g)$, so $\lambda(g)/k$
is the average information along the $k$ eigenvector directions. The
multivariate efficiency is $\mbox{Eff}^{\Lambda}(g)=I^{-1/2}\Lambda(g)I^{-1/2}$,
bounded above by the identity, with the score attaining the bound.

For the score-map image $g_{\hat{\mu}}$, the marginal score of the
statistic $\hat{\mu}$, the identity (\ref{eq:Lambda_score_map}) takes
the multivariate form $\Lambda(g_{\hat{\mu}})=I_{\hat{\mu}}$, the
Fisher information matrix of the marginal model of $\hat{\mu}$: writing
$s=g_{\hat{\mu}}$, $C=\mbox{Cov}(s,\nabla\ell)
=E\{s\,E(\nabla\ell\mid\hat{\mu})^{{\tt t}}\}=E(s\,s^{{\tt t}})=I_{\hat{\mu}}$
and $V(s)=I_{\hat{\mu}}$, so
$\Lambda=I_{\hat{\mu}}I_{\hat{\mu}}^{-1}I_{\hat{\mu}}=I_{\hat{\mu}}$.
The efficiency $\mbox{Eff}^{\Lambda}=I^{-1/2}I_{\hat{\mu}}I^{-1/2}$
equals the identity exactly when $\hat{\mu}$ is sufficient.

\section{The James--Stein example revisited\label{sec:JS-return}}

We now return to the example that motivated the paper. Sections \ref{sec:GE-vs-PE}
and \ref{sec:Lambda-vs-MSE} compared generalized estimators with point
estimators and $\Lambda$-information with MSE on general grounds, without
reference to James--Stein. We apply those comparisons here. The James--Stein
estimator was called paradoxical because it dominates maximum likelihood
in MSE while maximum likelihood is superior by other criteria; we show
that the paradox is an artifact of the assessment criterion, and that
under $\Lambda$-information it does not arise.

We consider $k$ independent observations from normal distributions with
unit variance and possibly different means:
\begin{equation}
M_{k}=\left\{ m(y)=(2\pi)^{-k/2}\exp(-\tfrac{1}{2}||y-\mu||^{2}),\ \mu\in\mathbb{R}^{k}\right\} .\label{eq:M}
\end{equation}
The ML estimate for $\mu$ is $\hat{\mu}=y\in\mathbb{R}^{k}$. The JS
estimate is
\[
\tilde{\mu}=\left(1-\frac{k-2}{||y||^{2}}\right)y,
\]
a shrinkage estimator that pulls the ML estimate toward the origin,
with the shrinkage factor depending on the squared norm of the observations.

For simplicity, we explore estimator properties using a subfamily.
Let $\mathbf{1}$ be the vector of ones in $\mathbb{R}^{k}$ and let
$\mu:M_{k}\rightarrow\mathbb{R}^{k}$ be the mean parameterization,
so
\[
M_{\mathbf{1}}=\left\{ m\in M_{k}:\mu(m)=\theta\mathbf{1},\theta\in\mathbb{R}\right\} 
\]
is a one-dimensional subfamily parameterized by $\theta$. This restriction
does not disadvantage JS; the MSE advantage of JS estimators increases
when means are similar. We use $\mu$ to indicate both a function on
$M_{k}$ and a point in its range. Our choice of $\theta$ and $k$
are motivated by \citet{Taketomi2013}'s gastric cancer data, where
the mean of $k=14$ estimates is $-0.245$ with mean standard error
$0.196$. Using unit standard error, we consider values near $\theta=0.245/0.196=1.25$.
Simulated quantities reported below are based on $1{,}000{,}000$
samples; quantities that follow from sufficiency are exact.

\subsection{Resolution of the paradox}

Assessed by MSE, the James--Stein estimator dominates: Table \ref{tab:MSE-for-James-Stein}
shows its MSE below that of maximum likelihood at every value of $\theta$,
the advantage greatest for $\theta$ near zero.

\begin{table}
\begin{tabular}{|c|r|r@{.}lrrr|}
\hline 
Estimator \textbackslash{} $\theta$:  & 0  & 0&50  & 1.25 & 2.00 & 2.50\tabularnewline
\hline 
\hline 
JS  & \textbf{2.00}  & \textbf{4}&\textbf{45}  & \textbf{9.58} & \textbf{11.83} & \textbf{12.53}\tabularnewline
\hline 
ML  & 14.00  & 14&00  & 14.00 & 14.00 & 14.00\tabularnewline
\hline 
\end{tabular}\caption{MSE for James--Stein (JS) and maximum likelihood (ML) estimators for
$k=14$ means having identical expectations $\theta\in\left\{ 0,0.5,1.25,2,2.5\right\} $.}\label{tab:MSE-for-James-Stein}
\end{table}

Assessed by information, the comparison reverses, and for maximum
likelihood no simulation is required. The ML estimator
$\hat{\mu}(y)=y$ is the identity statistic, trivially sufficient, so
under the score map its generalized estimator is the full-data score
$y-\mu$ and, by the multivariate form of (\ref{eq:Lambda_score_map})
recorded in Section \ref{subsec:bound},
$\Lambda(g_{\hat{\mu}})=I_{\hat{\mu}}=I_{k\times k}$ exactly, at every
distribution of $M_{k}$: scalar information $\lambda=14$ and mean
efficiency $1$ throughout. (The expectation map gives the same object here:
$g_{E(\hat{\mu})}=\hat{\mu}-E(\hat{\mu})=y-\mu$ as well.)

For James--Stein the two maps differ, and the score map exposes
structure that expectation-centering does not. Write $r=||y||$,
$u=y/r$, and $\rho(r)=r-(k-2)/r$, so that
$\tilde{\mu}(y)=\rho(r)\,u$. On the outer region $r>\sqrt{k-2}$,
$\rho$ is positive and increases from $0$ to $\infty$; on the inner
region $r<\sqrt{k-2}$, $\rho$ is negative---the negative-shrinkage
samples---and $|\rho|$ decreases from $\infty$ to $0$. Every nonzero
value $t$ of the JS statistic therefore has exactly two preimages,
\[
y_{1}=r_{1}\,\frac{t}{||t||},\qquad y_{2}=-\,r_{2}\,\frac{t}{||t||},
\qquad r_{2}<\sqrt{k-2}<r_{1},
\]
one on each branch, pointing in opposite directions. The JS statistic
is two-to-one: it collapses each such pair, discarding the binary
digit $B=\mathbf{1}\{||y||>\sqrt{k-2}\}$ that records which preimage
occurred. Its $\sigma$-field is that of the data minus one bit. The
collapse is visible at sample level as the familiar negative-shrinkage
anomaly: an inner-branch sample yields a JS estimate pointing opposite
to the data, and no user of $\tilde{\mu}$ alone can distinguish that
sample from its antipodal outer-branch partner reporting the same
value.

The information cost of the collapse is exact. The score of the
marginal model of $\tilde{\mu}$ is the conditional expectation of the
full score given the statistic,
\begin{equation}
s(t)=E\left(y-\mu\mid\tilde{\mu}=t\right)
=\pi(t)\,(y_{1}-\mu)+\{1-\pi(t)\}\,(y_{2}-\mu),\label{eq:JS_marginal_score}
\end{equation}
a two-atom average whose weight $\pi(t)$ is the conditional
probability of the outer branch, determined by the density of $y$ and
the Jacobian of the JS map at the two
preimages.\footnote{Both preimages are available in closed form: with
$c=||t||$, the radii solve $r^{2}\mp cr-(k-2)=0$, giving
$r_{1}=\{c+\sqrt{c^{2}+4(k-2)}\}/2$ and
$r_{2}=\{-c+\sqrt{c^{2}+4(k-2)}\}/2$, so that $r_{1}r_{2}=k-2$: the
preimages are inversions of one another through the sphere
$||y||^{2}=k-2$, with direction reversed. The branch weight is
$\pi(t)=w_{1}/(w_{1}+w_{2})$ with
$w_{i}=m(y_{i})/\{\rho'(r_{i})(c/r_{i})^{k-1}\}$,
$\rho'(r)=1+(k-2)/r^{2}$, and $m$ the density in (\ref{eq:M}). These
expressions make (\ref{eq:JS_marginal_score}) computable in closed
form at every sample, so $I_{\tilde{\mu}}=E(s\,s^{{\tt t}})$ is a
routine Monte Carlo average.} The variance decomposition of the full
score then gives
\begin{equation}
\Lambda(g_{\tilde{\mu}})=I_{\tilde{\mu}}
=I_{k\times k}-E\left\{ V\left(y\mid\tilde{\mu}\right)\right\},\label{eq:JS_info_loss}
\end{equation}
the full information minus the information destroyed by not knowing
the branch. Since $0<\pi(t)<1$ for every $t$, the deficit term in
(\ref{eq:JS_info_loss}) is nonzero and the JS efficiency is strictly
below one at every distribution in $M_{k}$. Where the deficit falls is
governed by the branch probabilities: when $||\mu||$ is large the
outer branch dominates, $\pi\approx1$, $s(t)\approx y-\mu$, and JS
retains nearly full information; near the shrinkage target $\mu=0$ the
event $||y||^{2}<k-2$ has appreciable probability, the two preimages
genuinely compete, and information is destroyed---precisely where the
MSE advantage of JS in Table \ref{tab:MSE-for-James-Stein} is
greatest. Table \ref{tab:Lambda-for-James-Stein} reports the
information, mean efficiency, and information deficit $k-\lambda$ for
JS along $M_{\mathbf{1}}$: the information and efficiency are
simulated from (\ref{eq:JS_marginal_score}), the deficit is computed
by numerical integration of the deficit term in
(\ref{eq:JS_info_loss}), and the ML rows are exact.

\begin{table}
\begin{tabular}{|c|c|c|c|c|c|}
\hline
Estimator & 0 & 0.50 & 1.25 & 2.00 & 2.50\tabularnewline
\hline
\hline
JS ($\lambda$) & 2.31 & 12.88 & 14.00 & 14.00 & 14.00\tabularnewline
\hline
ML ($\lambda$, exact) & \textbf{14} & \textbf{14} & \textbf{14} & \textbf{14} & \textbf{14}\tabularnewline
\hline
JS ($\lambda$-eff$/k$) & 0.17 & 0.92 & 1.00 & 1.00 & 1.00\tabularnewline
\hline
ML ($\lambda$-eff$/k$, exact) & \textbf{1} & \textbf{1} & \textbf{1} & \textbf{1} & \textbf{1}\tabularnewline
\hline
JS deficit ($k-\lambda$) & 11.69 & 1.12 & $5.5\times10^{-5}$ & $1.4\times10^{-12}$ & $1.6\times10^{-19}$\tabularnewline
\hline
\end{tabular}\caption{Scalar information $\lambda=\mbox{tr}\,\Lambda$ under the score map,
mean information efficiency, and information deficit $k-\lambda$ for
James--Stein (JS) and maximum likelihood (ML) estimators for $k=14$
means having identical expectations $\theta\in\left\{ 0,0.5,1.25,2,2.5\right\} $.
ML values are exact by sufficiency. JS information and efficiency are
computed by simulation from (\ref{eq:JS_marginal_score}); the deficit
row is computed by numerical integration of $E\{V(y\mid\tilde{\mu})\}$
in (\ref{eq:JS_info_loss}) and agrees with the simulation where the
simulation resolves it ($\theta\le1.25$), while for $\theta\ge2$ the
deficit lies far below simulation resolution.\label{tab:Lambda-for-James-Stein}}
\end{table}

The geography of the loss is stark. At the shrinkage target the
collapse destroys $11.7$ of the $14$ units of Fisher
information---efficiency $0.17$---while by $\theta=1.25$ the deficit
has fallen to $5.5\times10^{-5}$ and by $\theta=2.5$ to $10^{-19}$.
The JS statistic is inefficient everywhere in principle and nearly
sufficient everywhere in practice, except near its shrinkage target,
where the information cost is severe.

This is the dissolution of the paradox, and it is sharper than a mere
reversal of rankings. The James--Stein estimator purchases its MSE
advantage by discarding data: it is a two-to-one statistic, and the
discarded bit carries information about $\mu$ exactly where the
shrinkage is aggressive. MSE---a pointwise criterion, blind by the
diagnosis of Lemma \ref{lem:pointwise} to everything but the sampling
distribution at each member of $M$ separately---rewards the coarsening.
Information assessment sees it directly: the efficiency credited to
maximum likelihood is the efficiency of a sufficient statistic, and
the inefficiency charged to James--Stein is the exact information cost
of the collapse. There is no tension to resolve and no need to reconsider maximum likelihood.

\subsection{Estimation as the first step of inference\label{subsec:first-step}}

A natural objection is that the James--Stein estimator was designed
to minimize risk rather than to discriminate among the distributions
of $M_{k}$, so that an assessment based on discrimination holds it to
a standard it was not intended to meet. The objection mistakes what is
being assessed and why.
The question raised by the MSE comparison---Efron's question---is
whether the classical theory built on maximum likelihood needs
rethinking, and testing is part of that theory. A point estimate is
not, by itself, an inference. Classical practice is a two-step
procedure: first find a point estimator with good pointwise
properties, then build on it---standard errors, test statistics,
confidence intervals---to obtain inference procedures. An assessment
of the first step that ignores the demands of the second endorses
estimators that cannot support what is to be built on them.
Information measures exactly the resource that the second step
inherits: $I_{\tilde{\mu}}$ is the Fisher information the JS statistic
retains for every subsequent test and interval, however constructed,
and Table \ref{tab:Lambda-for-James-Stein} shows that resource
arriving largely destroyed precisely where the first-step advantage
lies. Generalized estimation dispenses with the two-step structure
altogether: familywise properties are demanded from the start, and a
generalized estimator is already a continuum of test statistics, one
for each distribution in $M$ (Section \ref{subsec:unification}).
Estimation and testing are not successive steps but a single object
read two ways.

\section{Conclusion}

For over six decades, the James--Stein estimator has stood as a challenge
to classical statistical theory, appearing to demonstrate that we
can systematically improve upon maximum likelihood by shrinking estimates
toward arbitrary points. This result has influenced countless applications,
from empirical Bayes methods to regularization in high-dimensional
statistics. See, for example, \citet{EfronMorris1973} and \citet{Tibshirani1996}.
Our analysis suggests these applications might benefit from reconsideration:
what appears as improvement through the lens of MSE may involve trade-offs
that become apparent only when viewed through the broader framework
of information and generalized estimation.

The framework introduced in \citet{Vos2024_2022} and extended by
\citet{VosWu2025} and \citet{Vos2025geometry} provides both a diagnosis
and a solution. The diagnosis:
MSE fails because it is defined to assess point estimators. Estimation,
as \citet[p.~250]{cox_theoretical_2000} note, is more
than this. MKL and risks based on other loss functions fail for the
same reason (Lemma \ref{lem:pointwise}). The solution: use generalized estimators to capture
aspects of estimation that are lacking in point estimators and use
$\Lambda$-information to assess these estimators. The maximum likelihood
estimator, viewed as a point estimate, is then replaced by its associated
generalized estimator---in the exponential family setting of the
James--Stein problem, the full-data score. So viewed, it attains the
Fisher information bound at every distribution.

Two centuries of statistical theory based on maximum likelihood need
not be rethought.  When we assess estimators beyond the numerical
value provided by their point estimate, the classical theory's wisdom
becomes apparent once again. The James--Stein ``paradox'' dissolves
when we stop using the wrong measuring stick.

\bibliographystyle{plainnat}
\bibliography{vos}

\end{document}